\documentclass[10pt,leqno]{article}
\usepackage{amssymb,amsbsy,amsmath,amsfonts,amssymb,amscd, mathrsfs}

\usepackage[english]{babel}
\usepackage[T1]{fontenc}
\usepackage{indentfirst}

\makeatletter
\@addtoreset{equation}{section}
\makeatother

\newtheorem{statement}{}[section]
\newtheorem{theoreme}[statement]{Theorem}
\newtheorem{lemme}[statement]{Lemma}

\newtheorem{proposition}[statement]{Proposition}
\newtheorem{definition}[statement]{Definition}

\newcommand\C{\mathbb C}

\newcommand\R{\mathbb R}
\newcommand\T{\mathbb T}
\newcommand\D{\mathbb D}

\renewcommand\H{\mathbb H} 
\newcommand\e{{\rm e}}
\renewcommand\P{\mathbb P}

\newcommand\ind{{\rm 1\kern-.30em I}}
\newcommand\qed{\hfill $\square$}
\renewcommand \Re{{\mathfrak R}{\rm e}\,}

\title{\bf  Some new properties of composition operators associated with lens maps}
\author{\it Pascal Lef\`evre,\it Daniel Li, \\
\it Herv\'e Queff\'elec, Luis Rodr{\'\i}guez-Piazza}

\date{\footnotesize \today}

\begin{document}

\maketitle

\noindent{\bf Abstract.} \emph{We give examples of results on composition operators connected with lens maps. The first two concern the 
approximation numbers of those operators acting on the usual Hardy space $H^2$.  The last ones are connected with Hardy-Orlicz and 
Bergman-Orlicz spaces $H^\psi$ and $B^\psi$, and provide a negative answer to the question of knowing if all composition operators which are weakly compact 
on a non-reflexive space are norm-compact.} 
\bigskip

\medskip



\section{Introduction}
 
We first recall the context of this work, which appears as a continuation of \cite{JFA}, \cite{JMAA}, \cite{MEMOIRS}, \cite{TRANS} and 
\cite{LIQUEROD}. \par

Let $\D$ be the open unit disk of the complex plane and ${\cal H} (\D)$ be the space of holomorphic functions on $\D$. To every analytic self-map 
$\varphi \colon \D \to \D$ (also called \emph{Schur function}), a linear map $C_\varphi \colon {\cal H} (\D) \to {\cal H} (\D)$ can be associated by 
$C_\varphi (f) = f \circ \varphi$.  This map is called the \emph{composition operator of symbol $\varphi$}. A  basic fact of the theory (\cite{Shap-livre}, 
page~13, or \cite{DUREN-livre}, Theorem~1.7) is \emph{Littlewood's subordination principle} which allows to show that every composition operator 
induces a \emph{bounded} linear map from the Hardy space $H^p$ into itself, $1 \leq p < \infty$. \par
\smallskip

In this work, we are specifically interested in a one-parameter family (a semi-group) of Schur functions: \emph{lens maps} $\varphi_\theta$, $0 < \theta < 1$, 
whose definition is given below. They turn out to be  very useful in the general theory of composition operators because they provide non trivial examples (for 
example, they generate compact and even Hilbert-Schmidt operators on the Hardy space $H^2$ \cite{Shap-livre}, page~27). The aim of this work is to illustrate 
that fact by new examples. \par
We show in Section~\ref{approx numbers} that, as operators on $H^2$, the approximation numbers of $C_{\varphi_\theta}$ behave as 
$\e^{-c_\theta \sqrt n}$. In particular, the composition operator $C_{\varphi_\theta}$ is in all Schatten classes $S_p$, $p > 0$. In Section~\ref{spreading}, 
we show that, when one ``spreads'' these lens maps, their approximation numbers become greater, and the associated composition operator 
$C_{\tilde \varphi_\theta}$ is in $S_p$ if and only if $p > 2 \theta$. In Section~\ref{counterexamples}, we answer to the negative a question of H.-O. Tylli: 
is it true that every weakly compact composition operator on a non-reflexive Banach function space is actually compact ? We show that there are composition 
operators on a (non-reflexive) Hardy-Orlicz spaces, which are weakly compact and Dunford-Pettis, though not compact and that there are composition operators 
on a non-reflexive Bergman-Orlicz space which are weakly compact but not compact. We also show that there are composition operators on a non-reflexive 
Hardy-Orlicz space which are weakly compact but not Dunford-Pettis. \par
\medskip
 
We give now the definition of lens maps (see \cite{Shap-livre}, page~27). 
\begin{definition} [Lens maps]
The \emph{lens map $\varphi_\theta \colon \D \to \D$ with parameter $\theta$}, $0 < \theta < 1$, is defined by: 
\begin{equation} \label{definition lens}
\qquad \varphi_{\theta} (z) = \frac{(1 + z)^\theta - (1 - z)^\theta }{(1 + z)^\theta + (1 - z)^\theta} \, \raise 1pt \hbox{,} \qquad z \in \D.
\end{equation}
\end{definition}

In a more explicit way, $\varphi_\theta$ is defined as follows.  Let $\H$ be the open right half-plane, and $T \colon \D \to \H$ be the (involutive) conformal 
mapping given by  
\begin{equation} \label{def de T}
T (z) =\frac{1 - z}{1 + z} \,\cdot
\end{equation} 
We denote by $\gamma_\theta$ the self-map of $\H$ defined by 
\begin{equation} \label{def de gamma} 
\gamma_{\theta} (w) = w^\theta = \e^{\theta \log w} ,
\end{equation} 
where $\log$ is the principal value of the logarithm and finally $\varphi_\theta \colon \D \to \D$ is defined by
\begin{equation} \label{lens} 
\varphi_{\theta} = T^{-1} \circ \gamma_\theta \circ T. 
\end{equation} 

Those lens maps form a continuous curve of analytic self-maps from $\D$ into itself, and an abelian semi-group for the composition of maps since we obviously 
have from (\ref{lens}) and the rules on powers that $\varphi_{\theta} (0) = 0$ and: 
\begin{equation}\label{semi} 
\varphi_{\theta} \circ \varphi_{\theta'} = \varphi_{\theta'} \circ \varphi_{\theta} = \varphi_{\theta\theta'} . 
\end{equation}

\noindent{\bf Acknowledgement.} The fourth-named author is partially supported by a Spanish research project MTM 2009-08934. \par
We thank the referee for a very quick answer and a very careful reading of this paper.
\goodbreak

\section{Approximation numbers of lens maps} \label{approx numbers}

For every operator $A \colon H^2 \to H^2$, we denote by 
\begin{displaymath} 
\qquad a_{n} (A) = \inf_{{\rm\ rank}\, R < n} \Vert A - R \Vert \,, \quad n = 1, 2, \ldots  
\end{displaymath} 
its \emph{$n$-th approximation number}. We refer to \cite{CA-ST-livre} for more details on those approximation numbers. \par  
\smallskip

Recall (\cite{Zhu-livre}, page~18) that the Schatten class $S_p$ on $H^2$ is defined by
\begin{displaymath} 
\qquad S_p = \{A \colon H^2 \to H^2 \, ; \  (a_{n}(A))_n \in \ell^p\}, \quad p > 0. 
\end{displaymath} 
$S_2$ is the Hilbert-Schmidt class and the quantity $\Vert A \Vert_p = \big(\sum_{n=1}^\infty (a_{n}(A))^p \big)^{1/p}$ is a Banach norm on $S_p$ for 
$p \geq 1$. \par
\smallskip

We can now state the following theorem:
\begin{theoreme}\label{approximation numbers} 
Let $0 < \theta < 1$ and $\varphi_\theta$ be the lens map defined in \eqref{definition lens}. There are positive constants $a, b, a', b'$ depending only on 
$\theta$ such that 
\begin{equation}\label{estimate} 
a'\,\e^{-b'\sqrt n} \leq a_{n} (C_{\varphi_\theta}) \leq a \, \e^{-b\sqrt n}.
\end{equation}
In particular, $C_{\varphi_\theta}$ lies in all Schatten classes $S_p$, $p > 0$.
\end{theoreme}

The lower bound in \eqref{approximation numbers} was proved in \cite{LIQUEROD}, Proposition~6.3. The fact that $C_{\varphi_\theta}$ lies in all Schatten 
classes was first proved in \cite{SHTA} under a qualitative form (see the very end of that paper). \par
The upper bound will be obtained below as a consequence of a result of O.~G.~Parfenov (\cite{Parf}). However, an idea of infinite divisibility, which may be 
used in other contexts, leads to a simpler proof, though it gives a worse estimate in \eqref{approximation numbers}: $\sqrt n$ is replaced by $n^{1/3}$. 
We shall begin by giving this proof, because it is quite short. It relies on the semi-group property \eqref{semi} and on an estimate of the Hilbert-Schmidt norm 
$\Vert C_{\varphi_\alpha}\Vert_2$  in terms of $\alpha$, as follows:
\begin{lemme} 
There exist numerical constants $K_1, K_2$ such that: 
\begin{equation} \label{Hilbert-Schmidt} 
\qquad \qquad \quad \frac{K_1}{1 - \alpha} \leq \Vert C_{\varphi_\alpha} \Vert_2 \leq \frac{K_2}{1 - \alpha} \, \raise 1pt \hbox{,} 
\qquad \quad \text{for all } 0 < \alpha < 1. \quad 
\end{equation}
In particular, we have 
\begin{equation} \label{weak} 
a_{n} (C_{\varphi_\alpha}) \leq \frac{K_2} {\sqrt n (1-\alpha)} \,\cdot 
\end{equation}
\end{lemme}

\noindent{\bf Proof.} The relation (\ref{weak}) is an obvious consequence of (\ref{Hilbert-Schmidt}) since
\begin{displaymath} 
n\big[a_{n} (C_{\varphi_\alpha}) \big]^2 
\leq \sum_{j=1}^n \big[a_{j}(C_{\varphi_\alpha})\big]^2 
\leq \sum_{j=1}^\infty \big[a_{j}(C_{\varphi_\alpha})\big]^2 
= \Vert C_{\varphi_\alpha}\Vert_2^{2} 
\leq \frac{K_{2}^2}{(1-\alpha)^2} \,\cdot
\end{displaymath} 

For the first part,  let $a = \cos (\alpha \pi/2) = \sin ((1-\alpha)\pi/2) \geq 1 - \alpha$ and let $\sigma = T (m)$ ($m$ is the normalized Lebesgue measure 
$dm (t) = dt/2 \pi$ on the unit circle) be the probability measure carried by the imaginary axis which satisfies:
\begin{displaymath} 
\int_{\H} f \, d\sigma = \int_{-\infty}^\infty f (i y) \frac{dy}{\pi (1 + y^2)} \,\cdot 
\end{displaymath} 
By definition, $T$, defined in \eqref{def de T}, is a unitary operator from $H^{2}(\D, m)$ into $H^{2}(\H, \sigma)$, and we easily obtain, setting 
$\gamma (y) =\gamma_{\alpha} (i y) = \e^{i (\pi/2) \alpha \, {\rm sign}\, (y)} \vert y \vert^\alpha$ (where ${\rm sign}$ is the sign of $y$ and 
$\gamma_\alpha$ is defined in \eqref{def de gamma}), that (see \cite{Shap-livre}, section~2.3):
\begin{align*}
\Vert C_{\varphi_\alpha} \Vert_2^{2} 
& = \int_\T \frac{dm}{1 - |\varphi_\alpha|^2} 
= \int_{\H} \frac{d\sigma}{1 - \big\vert \frac{1 - \gamma}{1 + \gamma} \big\vert^2} 
= \int_{\H} \frac{\vert 1 + \gamma \vert^2}{4 \, \Re \gamma} \,d\sigma \\
& =\int_{-\infty}^{+ \infty} \frac{\vert 1 + \gamma (y) \vert^2}{4 a \,\vert y \vert^\alpha} \,\frac{dy}{\pi (1+y^2)} \\
& \leq \frac{K}{1 - \alpha} \int_{0}^{+\infty} \frac{1 + y^{2\alpha}}{y^\alpha} \, \frac{dy}{1+y^2}
=\frac{2 K}{1 - \alpha} \int_{0}^{+\infty} \frac{y^\alpha}{1+y^2} \, dy \\
& \leq\frac{4 K}{(1 - \alpha)^2} \,\raise 1pt \hbox{,} 
\end{align*}
where $K$ is a numerical constant. This gives the upper bound in \eqref{Hilbert-Schmidt} and the lower one is obtained similarly.\hfill $\square$
\medskip

We can now finish the first proof of Theorem~\ref{approximation numbers}. Let $k$ be a positive integer and let 
\begin{displaymath} 
\qquad \alpha_k = \theta^{1/k},
\end{displaymath} 
so that $\alpha_{k}^k = \theta$. \par
Now use the well-known sub-multiplicativity $a_{p + q - 1} (vu) \leq a_{p}(v) \,a_{q}(u)$ of approximation numbers (\cite{Pis-livre}, page 61), as well as 
the semi-group property \eqref{semi} (which implies $C_{\varphi_\theta} = C_{\varphi_{\alpha_k}}^{k}$), and \eqref{weak}. We see that: 
\begin{displaymath} 
a_{kn} (C_{\varphi_\theta}) = a_{kn} (C_{\varphi_{\alpha_k}}^k) \leq \big[a_{n} (C_{\varphi_{\alpha_k}})\big]^k 
\leq  \bigg[\frac{K_2}{(1 - \alpha_k) \sqrt n}\bigg]^k.  
\end{displaymath} 
Observe that 
\begin{displaymath} 
1 - \alpha_k \geq \frac{1 - \alpha_{k}^k}{ k} = \frac{1 - \theta}{k} \,\cdot
\end{displaymath} 
We then get, $c = c_\theta$ denoting a constant which only depends on $\theta$:
\begin{displaymath} 
a_{kn} (C_{\varphi_\theta}) \leq \Big(\frac{k}{c\sqrt n}\Big)^k. 
\end{displaymath} 
Set $d = c/ \e$ and take $k = d \sqrt n$, ignoring the questions of integer part. We obtain:
\begin{displaymath} 
a_{d n^{3/2}} (C_{\varphi_\theta}) \leq \e^{- k} = \e^{- d \sqrt n}. 
\end{displaymath} 
Setting $N = d n^{3/2}$, we get 
\begin{equation} \label{mieux}
a_{N} (C_{\varphi_\theta}) \leq a\, \e^{- b N^{1/3}} 
\end{equation} 
for an appropriate value of $a$ and $b$ and for any integer $N\geq 1$. 
This ends our first proof, with an exponent slightly smaller that the right one ($1/3$ instead  of $1/2$), yet more than sufficient to prove that 
$C_{\varphi_\theta} \in \cap_{p > 0} S_p$. \hfill $\square$ 
\medskip

\noindent{\bf Remark.} Since the estimate \eqref{weak} is rather crude, it might be expected that, using \eqref{mieux}, and iterating the process, we could obtain 
a better one. This is not the case, and this iteration leads to \eqref{mieux} and the exponent $1/3$ again (with different constants $a$ and $b$). \par
\medskip

\noindent{\bf Proof of Theorem~\ref{approximation numbers}.} This proof will give the correct exponent $1/2$ in the upper bound. Moreover, it works 
more generally for Schur functions whose image lies in polygon inscribed in the unit disk. This upper bound appears, in a different context and under a very 
cryptic form, in \cite{Parf}.  First note the following simple lemma.
\begin{lemme}\label{Simple} 
Suppose that $a,b \in\D$ satisfy $\vert a - b \vert \leq M \min (1 - \vert a \vert, 1 - \vert b \vert)$, where $M$ is a constant. Then:
\begin{displaymath} 
d (a,b) \leq \frac{M}{\sqrt{M^2 +1}} := \chi < 1 .
\end{displaymath} 
\end{lemme} 
Here $d$ is the pseudo-hyperbolic distance defined by:
\begin{displaymath} 
\qquad d (a, b) = \Big\vert \frac {a - b}{1 - \overline{a} b} \Big\vert \,\quad a, b \in \D. 
\end{displaymath} 

\noindent{\bf Proof.} Set $\delta = \min (1 - \vert a \vert, 1 - \vert b \vert)$. We have the  identity
\begin{displaymath} 
\frac{1}{d^{2} (a,b)} - 1 = \frac{(1 - \vert a \vert^2) (1 - \vert b \vert^2)}{\vert a - b \vert^2} 
\geq \frac{(1 - \vert a \vert) (1 - \vert b \vert)}{\vert a - b \vert^2} 
\geq \frac{\delta^2}{M^{2} \delta^2} = \frac{1}{M^2} \, \raise 1pt \hbox{,}
\end{displaymath} 
hence the lemma. \hfill $\square$
\medskip

The second lemma gives an upper bound for $ a_{N} (C_\varphi)$. In this lemma, $\kappa$ is a numerical constant, $S(\xi, h)$ the usual 
\emph{pseudo-Carleson window} centered at $\xi \in \T$ (where $\T = \partial \D$ is the unit circle) and of radius $h$ ($0 < h < 1$), defined by:
\begin{equation} \label{pseudo-window} 
S (\xi, h) = \{z \in \D \, ; \ \vert z - \xi \vert \leq h\},
\end{equation} 
and $m_\varphi$ is the pull-back measure of $m$, the normalized Lebesgue measure on $\T$, by $\varphi^\ast$. Recall that if $f \in {\cal H} (\D)$, one sets  
$f_{r} (\e^{it}) = f (r \e^{it})$ for $0 < r < 1$ and, if the limit exists $m$-almost everywhere, one sets:
\begin{equation}\label{radial} 
f^*(e^{it})=\lim_{r\to 1^-}f (r \e^{it}).
\end{equation}   
Actually, we shall do write $f$ instead of $f^\ast$. 
Recall that a measure $\mu$ on $\overline \D$ is called a \emph{Carleson measure} if there is a constant $c > 0$ such that 
$\mu \big[ \overline {S (\xi, h)} \big] \leq c\, h$ for all $\xi \in \T$. Carleson's embedding theorem says that $\mu$ is a Carleson measure if and only if the 
inclusion map from $H^2$ into $L^2 (\mu)$ is bounded (see \cite{DUREN-livre}, Theorem~9.3, for example). 
\begin{lemme}\label{Embedding} 
Let $B$ be a Blaschke product with less than $N$ zeroes (each zero being counted with its multiplicity). Then, for every Schur function $\varphi$, one has:
\begin{equation} \label{Blaschke} 
a_{N}^2 := \big[ a_{N} (C_{\varphi}) \big]^2 
\leq \kappa^2 \sup_{0 < h < 1, \xi \in \T} \ \frac{1}{h} \int_{\overline{S (\xi, h)}} \vert B \vert^2 \, dm_{\varphi} ,
\end{equation}
for some universal constant $\kappa > 0$.
\end{lemme}
\noindent{\bf Proof.} The subspace $B H^2$ is of codimension $\leq N - 1$. Therefore, 
$a_N = c_{N} (C_\varphi) \leq \big\| {C_\varphi}_{\mid BH^2} \big\|$, where the $c_{N}$'s are the Gelfand numbers (see \cite{CA-ST-livre}), and 
where we used the equality $a_N = c_N$ occurring in the Hilbertian case (see \cite{CA-ST-livre}). Now, since $\Vert B f \Vert_{H^2} = \Vert f \Vert_{H^2}$ 
for any $f \in H^2$, we have:
\begin{align*} 
\big\| {C_\varphi }_{\mid B H^2} \big\|^2 
& =\sup_{\Vert f \Vert_{H^2} \leq 1} \int_\T \vert B \circ \varphi \vert^2 \, \vert f \circ \varphi \vert^2 \,dm 
=\sup_{\Vert f \Vert_{H^2} \leq 1} \int_{\overline \D} \vert B \vert^2\vert f \vert^2 \, dm_\varphi \\
& =\Vert R_\mu\Vert^2 , 
\end{align*} 
where $\mu = \vert B \vert^{2} m_\varphi$ and where $R_\mu \colon H^2 \to L^{2} (\mu)$ is the restriction map. Of course, $\mu$ is a Carleson measure 
for $H^2$ since $\mu \leq m_\varphi$. Now, Carleson's embedding theorem says us that 
$\Vert R_\mu \Vert^2 \leq \kappa^2 \sup_{0 < h < 1, \xi \in \T} \frac{\mu [\overline{S (\xi, h)}]}{h}$ (see \cite{DUREN-livre}, Remark after the proof of 
Theorem~9.3, at the top of page~163; actually, in that book, Carleson's windows $W (\xi, h)$ are used instead of pseudo-Carleson's windows $S (\xi, h)$, but that 
does not matter, since $W (\xi, h) \subseteq S (\xi, 2h)$: if $r \geq 1 - h$ and $|t - t_0| \leq h$, then 
$|r \e^{it} - \e^{i t_0}|\leq |r \e^{it} - \e^{it}| + |\e^{it} - \e^{it_0}| \leq 2h$). That ends the proof of Lemma~\ref{Embedding}. \hfill $\square$ 
\bigskip

The following lemma takes into account the behaviour of $\varphi_\theta (\e^{it})$, and will be useful to us in Section~\ref{spreading} as well. The notation 
$u (t) \approx v (t)$ means that $a \, u (t) \leq v (t) \leq b \, u (t)$, for some positive constants $a, b$.

\begin{lemme} \label{Behaviour}  
Set $\gamma (t) = \varphi_\theta (\e^{it}) = \vert \gamma (t) \vert \, \e^{iA (t)}$, with $- \pi \leq t \leq \pi$, and $- \pi \leq A(t) \leq \pi$. 
Then, for $0 \leq \vert t \vert , \vert t' \vert \leq \pi/2$, one has:
\begin{equation} \label{Lip}   
\vert 1 - \gamma (t) \vert \approx 1 - \vert \gamma (t) \vert \approx \vert t \vert^\theta \quad \text{and}\quad 
\vert \gamma (t) - \gamma (t') \vert \leq K\, \vert t - t' \vert^\theta .
\end{equation}
Moreover, we have for $\vert t \vert \leq \pi/2$: 
\begin{equation} \label{pointone} 
A (t) \approx \vert t \vert^\theta \quad \text{and} \quad A' (t) \approx \vert t \vert^{\theta - 1}. 
\end{equation}
\end{lemme}

\noindent {\bf Proof.} First, recall that 
\begin{displaymath} 
\varphi_\theta (z) = \frac{(1+ z)^\theta - (1 - z)^\theta}{(1 + z)^\theta + (1 - z)^\theta} \, \raise 1pt \hbox{,}
\end{displaymath} 
so that $\varphi_\theta (\overline{z}) = \overline{\varphi_\theta (z)}$ and $\varphi_\theta (- z) = - \varphi_\theta (z)$. It follows that 
$\gamma (- t) = \overline{\gamma (t)}$ and $\gamma (t + \pi) = - \gamma (t)$, so that we may assume $0 \leq t, t' \leq \pi/2$. Then, we have more precisely, 
setting $c = \e^{- i \theta \pi/2}$, $s = \sin (\theta \pi/2)$ and $\tau = \big( \tan (t/2) \big)^\theta$: 
\begin{displaymath} 
\gamma (t) = \frac{(\cos t/2)^\theta  - \e^{-i \theta\pi/2} (\sin t/2)^\theta} {(\cos t/2)^\theta + \e^{-i \theta \pi/2} (\sin t/2)^\theta} 
= \frac{1 - c \tau}{1 + c \tau} = \frac{1 - \tau^2}{\vert 1 + c \tau\vert^2} +\frac{2 i s \tau}{\vert 1 + c \tau\vert^2} \, \raise 1pt \hbox{,}
\end{displaymath} 
after a simple computation, since $(1 + \e^{it})^\theta = \e^{i t \theta/2} (2 \cos t/2)^\theta$ and 
$(1 - \e^{it})^\theta = \e^{- i \theta \pi/2} \, \e^{i t \theta/2} (2 \sin t/2)^\theta$. Note by the way that 
\begin{displaymath} 
\varphi_\theta (1) = 1 \,; \quad \varphi_\theta (i) = i \tan (\theta \pi /4) \, ; \quad 
\varphi_\theta (-1) = - 1 \, ; \quad \varphi_\theta (- i) = - i \tan (\theta \pi/4).
\end{displaymath} 

Now, observe that $2 \geq \vert 1 + c \tau \vert \geq \Re (1 + c \tau) \geq 1$ and therefore that 
\begin{displaymath} 
\vert 1 - \gamma (t) \vert = \bigg\vert \frac{2 c \tau}{1 + c \tau} \bigg\vert \approx \tau \approx t^\theta \,, 
\end{displaymath} 
and similarly for $1 - \vert \gamma (t) \vert$ since $1 - \vert \gamma (t) \vert^2 = \frac{4 (\Re c) \, \tau}{\vert 1 + c \tau \vert^2}$.
The  relation \eqref{Lip} clearly follows. To prove \eqref{pointone}, we just have to note that, for $0 \leq t \leq \pi/2$, we have  
$A (t) = \arctan \frac{2 s\tau}{1 - \tau^2} \,\cdot$ \hfill $\square$ 
 \bigskip

Now, we prove Theorem~\ref{approximation numbers} in the following form (in which $q = q_\theta$ denotes a positive constant smaller than one), 
which is clearly sufficient. 
 \begin{equation} \label{Form} 
a_{4 N^{2} + 1} \leq K q^N. 
\end{equation}
The proof will come from an adequate choice of a Blaschke product of length $4 N^{2}$, with zeroes on the curve 
$\gamma (t) = \varphi_\theta (\e^{it})$, $-\pi \leq t \leq \pi$. Let $t_k = \pi 2^{-k}$ and $p_k = \gamma (t_k)$, with $1 \leq k \leq N$, so that the points 
$p_k$ are all in the first quadrant. We reflect them through the coordinate axes, setting: 
\begin{displaymath} 
\qquad q_k = \overline{p_k}, \quad r_k = - p_k, \quad s_k = - q_k, \qquad 1 \leq k \leq N. 
\end{displaymath} 
Let now $B$ be the Blaschke product having a zero of order $N$ at each of the points $p_k, q_k, r_k, s_k$, namely:
\begin{displaymath} 
B (z) = \prod_{k = 1}^N \bigg[ \frac{z - p_k}{1 - \overline{p_k} z} \cdot  
\frac{z - q_k}{1 - \overline{q_k} z} \cdot \frac{z - r_k}{1 - \overline{r_k} z} \cdot \frac{z - s_k}{1 - \overline{s_k} z} \bigg]^N.  
\end{displaymath} 
This Blaschke product  satisfies, by construction, the symmetry relations: 
\begin{equation} \label{symmetry} 
\quad B (\overline{z}) = \overline{B (z)} \,, \quad B (- z) = B (z). 
\end{equation}
 Of course, $\vert B\vert = 1$ on the boundary of $\D$, but $\vert B \vert$ is small on a large portion of the curve $\gamma$, as expressed by the following 
 lemma.
\begin{lemme} \label{Goodportion} 
For some constant $\chi = \chi_\theta < 1$, the following estimate holds:
\begin{equation} \label{Estimate} 
t_N \leq t \leq t_1 \quad \Longrightarrow \quad \vert B (\gamma (t) ) \vert \leq \chi^N. 
\end{equation}
\end{lemme}
\noindent{\bf Proof.} Let $t_N \leq t \leq t_1$ and $k$ such that $t_{k+1} \leq t \leq t_k$.  Let $B_k (z) =  \frac{z - p_k}{1 - \overline{p_k} z}$. Then, 
with help of Lemma~\ref{Behaviour}, we see that the assumptions of Lemma~\ref{Simple} are satisfied with $a = \gamma (t)$ and $b = \gamma(t_k)$, since 
$\vert t - t_k \vert \leq t_k - t_{k+1} =\pi 2^{- k - 1}$, so that $\min (1 - \vert a \vert, 1 - \vert b \vert) \approx t_{k}^\theta \approx 2^{- k\theta}$ and
hence, for some constant $M$: 
\begin{displaymath} 
\qquad
\vert a - b \vert \leq K \, \vert t - t_k \vert^\theta \leq K 2^{- k \theta} \leq M\, \min (1 - \vert a \vert, 1 - \vert b \vert) \,.
\end{displaymath} 

We therefore have, by definition, and by Lemma~\ref{Simple}, where we set $\chi = M/ \sqrt{M^2 + 1}$: 
\begin{displaymath} 
\vert B_{k} (\gamma (t) ) \vert =  d (\gamma(t), p_k) \leq \chi < 1.
\end{displaymath} 
It then follows from the definition of $B$ that:
\begin{displaymath} 
\vert B (\gamma (t) ) \vert \leq \vert B_{k} (\gamma(t) ) \vert^{N} \leq \chi^{N} ,
\end{displaymath} 
and that ends the proof of Lemma~\ref{Goodportion}. \hfill $\square$
\bigskip

Now fix $\xi \in \T$ and $0 < h \leq 1$. By interpolation, we may assume that $h = 2^{- n \theta}$. By symmetry, we may assume that $\Re \xi \geq 0$ and 
$\Re \gamma (t) \geq 0$, i.e. $\vert t \vert \leq \pi/2$. Then, since $\varphi_\theta (\D)$ is contained in the symmetric angular sector of vertex $1$ and opening 
$\theta \pi < \pi$, there is a constant $K > 0$ such that 
$\vert 1 - \gamma (t) \vert \leq K (1 - \vert \gamma (t) \vert)$. The only pseudo-windows $S(\xi, h)$ giving an integral not equal to zero in the estimation 
\eqref{Blaschke} of Lemma~\ref{Embedding} satisfy $\vert \xi - 1 \vert \leq (K + 1) h$. Indeed, suppose that $\vert \gamma (t) - \xi \vert \leq h$.  Then 
$1 - |\gamma (t) | \leq |\gamma (t) - \xi | \leq h$ and $\vert 1 - \gamma (t) \vert \leq K (1 - \vert \gamma (t) \vert) \leq  K h$. If  
$\vert \xi - 1 \vert > (K + 1) h$, we should have   $\vert \gamma (t) - \xi \vert \geq \vert \xi - 1 \vert - \vert \gamma (t) - 1 \vert > (K + 1) h - K h = h$, 
which is impossible.  Now, for such a window,  we have by definition of $m_\varphi$:
\begin{align*}
\int_{S (\xi, h)} \vert B \vert^2 \, dm_{\varphi_\theta} 
& = \int_{\vert \gamma (t) - \xi \vert \leq h} \vert B (\gamma (t) ) \vert^2 \, \frac{dt}{2 \pi}  
\leq \int_{\vert \gamma (t) - 1 \vert \leq (K + 2) h} \vert B (\gamma (t) ) \vert^2 \, \frac{dt}{2 \pi} \\ 
& \leq \int_{\vert t \vert \leq D t_n} \vert B (\gamma (t) ) \vert^2 \, \frac{dt}{2\pi} \mathop{=}^{def} I_h , 
\end{align*}
since $\vert \gamma (t) - 1 \vert \leq \vert \gamma (t) - \xi \vert + \vert \xi - 1 \vert \leq h + (K + 1) h$ and since 
$\vert \gamma (t) - 1 \vert \geq a \vert t \vert^\theta$ and $\vert \gamma (t) - 1 \vert \leq (K + 2) h$ together imply $\vert t \vert \leq D t_n$, where 
$D > 1$ is another constant (recall that $h = 2^{- n \theta} = (t_n / \pi)^\theta$).  \par 
\smallskip

To finish the discussion,  we separate two cases. 
\par

1) If $n \geq N$, we simply majorize $\vert B \vert$ by $1$. We set $q_1 = 2^{\theta - 1} < 1$ and get: 
\begin{displaymath} 
\frac{1}{h}\, I_h \leq \frac{1}{h} \int_{- D t_n}^{D t_n} \vert B (\gamma (t) ) \vert^2 \, \frac{dt}{2 \pi}  
\leq \frac{2 D t_n}{2 \pi h} = D q_{1}^n \leq D \, q_{1}^N.
\end{displaymath} 
\par

2) If $n \leq  N - 1$, we write:
\begin{displaymath} 
\frac{1}{h} \, I_h = \frac{2}{ h} \int_{0}^{D t_N} \vert B (\gamma (t) ) \vert^2 \, \frac{dt}{2 \pi} 
+ \frac{2}{h}\int_{D t_{N}}^{D t_n} \vert B (\gamma (t) ) \vert^2 \, \frac{dt}{2 \pi} := J_N + K_N.
\end{displaymath} 
The term $J_N$ is estimated above: $J_N \leq D \, q_{1}^N$. The term $K_N$ is estimated through Lemma~\ref{Goodportion}, which gives us:
\begin{displaymath} 
K_N \leq 2^{n \theta} \, \frac{2 D t_n}{2 \pi}\,  \chi^{2N}  \leq D \, \chi^{2N} , 
\end{displaymath} 
since $t_n 2^{n \theta} \leq \pi$, due to the fact that $\theta < 1$. \par
\smallskip

If we now apply Lemma~\ref{Embedding} with $q = \max (q_1, \chi^2)$ and with $N$ changed into $4 N^2 + 1$, we obtain \eqref{Form}, by changing 
the value of the constant $K$ once more. This ends the proof of Theorem~\ref{approximation numbers}. \hfill $\square$
\bigskip 

Theorem~\ref{approximation numbers} has the following consequence (as in \cite{Shap-livre}, page~29). 
\begin{proposition}
Let $\varphi$ be a \emph{univalent} Schur function and assume that $\varphi (\D)$ contains an angular sector centered on the unit circle and with opening 
$\theta \pi$, $0 < \theta < 1$. Then $a_n (C_\varphi) \geq a \, \e^{- b \sqrt n}$, $n = 1, 2, \ldots$, for some positive constants $a$ and $b$, depending only 
on $\theta$.
\end{proposition}
\noindent{\bf Proof.} We may assume that this angular sector is centered at $1$. By hypothesis, $\varphi (\D)$ contains the image of the ``reduced'' lens map 
defined by $\tilde \varphi_\theta (z) = \varphi_\theta ( (1 + z)/ 2)$. Since $\varphi$ is univalent, there is a Schur function $u$ such that 
$\tilde \varphi_\theta = \varphi \circ u $. Hence $C_{\tilde \varphi_\theta} = C_\varphi \circ C_u$ and 
$a_n (C_{\tilde \varphi_\theta}) \leq \| C_u\| \, a_n (C_\varphi)$. Theorem~\ref{approximation numbers} gives the result, since the calculations for 
$\tilde \varphi_\theta$ are exactly the same as for $\varphi_\theta$ (because they are equivalent as $z$ tends to $1$). \qed 
\medskip

The same is true if $\varphi$ is univalent and $\varphi (\D)$ contains a polygon with vertices on $\partial \D$. 


\section{Spreading the lens map} \label{spreading} 

In \cite{JFA}, we studied the effect of the multiplication of a Schur function $\varphi$ by the singular inner function 
$M (z) = \e^{- \frac{1 + z}{1 - z}}$, and observed that this multiplication spreads the values of the radial limits of the symbol and  lessens the maximal occupation 
time for Carleson windows. In some cases this improves the compactness or membership to Schatten classes of $C_\varphi$. More precisely, we proved 
the following result.  

\begin{theoreme} [\cite{JFA}, Theorem~4.2] 
For every $p > 2$, there exist two Schur functions $\varphi_{1}$ and $\varphi_2 = \varphi_{1} M$ such that 
$\vert \varphi_{1}^* \vert = \vert \varphi_{2}^* \vert$ and  $C_{\varphi_1} \colon H^2 \to H^2$ is not compact, but 
$C_{\varphi_2} \colon H^2 \to H^2$ is in the Schatten class $S_p$.
\end{theoreme}

Here, we will meet the \emph{opposite phenomenon}: the symbol $\varphi_1$ will have a fairly big associated maximal function 
$\rho_{\varphi_1}$, but will belong to all Schatten classes since it ``visits'' a bounded number of windows (meaning that there exists an integer $J$ such that, 
for fixed  $n$, at most $J$ of the $W_{n, j}$ are visited by $\varphi^{*} (\e^{it})$). The spread symbol will have an improved maximal function, but will visit 
all windows, so that its membership in Schatten classes will be degraded. More precisely, we will prove that
\begin{theoreme} \label{degrade} 
Fix $0 < \theta < 1$. Then there exist two Schur functions $\varphi_{1}$ and $\varphi_2$ such that: \par
\smallskip
1) $C_{\varphi_1} \colon H^2 \to H^2$ is in all Schatten classes $S_p$, $p > 0$, and even $a_{n} (C_{\varphi_1}) \leq a\, \e^{- b \sqrt n}$; \par
\smallskip 
2) $\vert \varphi_{1}^* \vert = \vert \varphi_{2}^* \vert$; \par
\smallskip
3) $C_{\varphi_2} \in S_p$ if and only if $p > 2 \theta$; \par
\smallskip
4) $a_{n} (C_{\varphi_2}) \leq K \, (\log n / n)^{1 /2 \theta}$, $n = 2, 3, \ldots$. 
\end{theoreme}

Of course, it would be better to have a good lower bound for $a_{n} (C_{\varphi_2})$, but we do not succeed in finding it yet.\par
\smallskip

\noindent{\bf Proof.} First observe that $C_{\varphi_1} \in S_2$, so that $C_{\varphi_2}\in S_2$ too, since 
$\vert \varphi_{1}^* \vert = \vert \varphi_{2}^* \vert$ and since the membership of $C_\varphi$  in $S_2$ only depends on the modulus of 
$\varphi^{*}$ because it amounts to (\cite{Shap-livre}, page 26): 
\begin{displaymath} 
\int_{-\pi}^\pi \frac{dt}{1 - \vert \varphi^{*} (\e^{it})\vert} < \infty. 
\end{displaymath} 
Theorem~\ref{degrade} says that we can hardly have more. We first prove a lemma. Recall (see \cite{JFA}, for example) that the 
\emph{maximal Carleson function} $\rho_\varphi$ of a Schur function $\varphi$ is defined, for $0 < h < 1$, by:
\begin{equation} 
\rho_\varphi (h) = \sup_{|\xi| = 1} m_\varphi[ S (\xi, h)] . 
\end{equation} 
\begin{lemme} \label{simplet} 
Let $0 < \theta < 1$. Then,  the maximal function $\rho_{\varphi_\theta}$ of $\varphi_\theta$ satisfies 
$\rho_{\varphi_\theta} (h) \leq K^{1/\theta} (1 - \theta)^{- 1/\theta} h^{1/\theta}$ and, moreover,
 \begin{equation} \label{Edensor case} 
\rho_{\varphi_\theta} (h) \approx h^{1/\theta} .
\end{equation} 
\end{lemme}
\noindent{\bf Proof of the lemma.} Let $0 < h < 1$ and $\gamma (t) = \varphi_\theta (\e^{it})$. $K$ and $\delta$ will denote constants which can change 
from a formula to another. We have, for $\vert t \vert \leq \pi/2$:
\begin{align*}
1 - \vert \gamma (t) \vert^2 
& = \frac{4 (\Re c) \tau}{\vert 1 + c \tau \vert^2} 
\geq \delta \cos (\theta \pi/2) \frac{\tau}{\vert 1 + c \tau \vert^2} 
\geq \delta (1 - \theta) \frac{\tau}{\vert 1 + c \tau \vert^2} \\
& \geq \delta (1 - \theta) \vert t \vert^\theta .
\end{align*}
Hence, we get, from Lemma~\ref{Behaviour}:
\begin{align*}
\rho_{\varphi_\theta} (h) 
& \leq 2 \, m (\{1 - \vert \gamma (t) \vert \leq h \text{ and } \vert t \vert \leq \pi/2\}) 
\leq  2 \, m (\{(1 - \theta) \delta \vert t \vert^\theta \leq K h\}) \\
& \leq K^{1/\theta} (1 - \theta)^{- 1/ \theta} h^{1/\theta}.
\end{align*}
Similarly, we have:
\begin{displaymath} 
\rho_{\varphi_\theta} (h) 
\geq m_{\varphi_\theta}[S(1, h)] 
\geq m (\{\vert 1 - \gamma (t) \vert \leq h\}) 
\geq m(\{\vert t\vert^\theta\leq K h\}) 
\geq K h^{1/\theta} ,
\end{displaymath} 
and that ends the proof of the lemma.  \hfill $\square$
\medskip

Going back to the proof of Theorem~\ref{degrade}, we take $\varphi_1 = \varphi_\theta$ and $\varphi_2 (z) = \varphi_{1} (z) M (z^2)$. We use 
$M(z^2)$ instead of $M (z)$ in order to treat the points $- 1$ and $1$ together. 
\par

The first two assertions are clear. 
For the third one, we define the dyadic Carleson windows, for $n = 1, 2, \ldots, j = 0, 1, \ldots, 2^{n} - 1$, by: 
\begin{displaymath} 
W_{n, j} = \{z \in \D \, ; \ 1 - 2^{- n} \leq \vert z \vert < 1 \text{ and } (2 j \pi) 2^{- n} \leq \arg (z) < (2 (j + 1)) \pi) 2^{- n} \} .
\end{displaymath} 
Recall (see \cite{JFA}, Proposition~3.3) the following proposition, which is a variant of Luecking's criterion (\cite{LUECKING}) for membership in a Schatten 
class, and which might also be used to give a third proof of the membership of $C_{\varphi_\theta}$ in all Schatten classes $S_p$, $p>0$, although the 
first proof turns out to be more elementary.
\begin{proposition} [\cite{LUECKING}, \cite{JFA}] \label{lulu} 
Let $\varphi$ be a Schur function and $p > 0$ a positive real number. Then $C_\varphi \in S_p$ if and only if 
\begin{displaymath} 
\sum_{n = 1}^\infty \sum_{j = 0}^{2^{n}-1} \big[2^n m_{\varphi} (W_{n, j}) \big]^{p/2} < \infty.
\end{displaymath} 
\end{proposition}

We apply this proposition with $\varphi = \varphi_2$, which satisfies, for $0 < \vert t \vert \leq \pi / 2$, the following relation: 
\begin{displaymath} 
\varphi (\e^{it}) = \vert \gamma (t) \vert \e^{i [A(t) - \cot ( t )]} \mathop{=}^{def} \vert \gamma (t) \vert \e^{i B (t)}, 
\end{displaymath} 
where $\gamma (t) = \varphi_{1} (\e^{i t})$ and (using Lemma~\ref{Behaviour}):
\begin{equation} \label{goodbehaviour} 
0 < \vert t \vert \leq \pi / 2 \ \Longrightarrow \ B (t) = \Gamma (t) - \frac{1}{t} \, \raise 1pt \hbox{,} 
\quad \text{ with } \Gamma (t) \approx \vert t \vert^\theta  \text{and } 
\Gamma ' (t) \approx \vert t \vert^{\theta - 1}.
\end{equation}
It clearly follows from \eqref{goodbehaviour} that the  function $B$ is increasing on some interval $[- \delta, 0[$ where $\delta$ is a positive numerical constant. 
Let us fix a positive integer $q_0$ such that $- \pi / 2 \leq t < 0 $ and
\begin{displaymath} 
B (t) \geq 2 q_0 \pi \quad \Longrightarrow \quad t \geq - \delta.
\end{displaymath} 
Fix a Carleson window $W_{n, j}$ and let us analyze the set $E_{n, j}$ of those $t$'s such that $\varphi (\e^{it})$ belongs to $W_{n, j}$. Recall that 
$m_{\varphi}(W_{n, j}) = m (E_{n, j})$. The membership in $E_{n, j}$ gives two constraints.\par
\smallskip

1) \emph{Modulus constraint}. We must have $\vert \gamma (t) \vert \geq 1 - 2^{- n}$, and therefore $\vert t \vert \leq K 2^{- n/\theta}$.
\par\smallskip
  
2) \emph{Argument constraint}. Let us set $\theta_{n, j} = (2j + 1) \pi 2^{- n}$, $h = \pi 2^{- n}$ and $I_{n, j} = (\theta_{n, j} - h, \theta_{n, j} + h)$. 
The angular constraint $\arg \varphi (\e^{it}) \in I_{n, j}$ will be satisfied if  $t < 0$ and
\begin{displaymath} 
B (t) \in \bigcup_{q \geq q_0} \big[\theta_{n, j} - h + 2 q \pi, \, \theta_{n, j} + h + 2 q \pi \big] := \bigcup_{q \geq q_0} J_{q} (h) := F.
\end{displaymath} 
We have $F \subset [2 q_0 \pi, \infty [$, and so $B (t) \in F$ and $t < 0$ imply $t \geq - \delta$. Set:
\begin{displaymath} 
E = \bigcup_{q \geq q_0} \big[B^{-1} (\theta_{n,j } - h + 2 q \pi), \, B^{-1} (\theta_{n, j} + h + 2 q \pi) \big] 
:= \bigcup_{q \geq q_0} I_{q} (h) \subset [- \delta, 0[.
\end{displaymath} 
The intervals $I_q$'s are disjoint, since $\theta_{n, j} + 2 (q+1) \pi - h > \theta_{n, j} + 2 q \pi + h$ and since $B$ increases on $[-\delta, 0[$. 
Moreover, $t \in E$ implies that $B (t) \in F$, which in turn implies that $\arg \varphi (\e^{it}) \in I_{n, j}$. Using Lemma~\ref{Behaviour}, we can find 
positive constants $c_1, c_2$ such that: 
\begin{displaymath} 
q \geq q_0 \quad \Longrightarrow \quad  - c_1/q \leq \min I_{q} (h) \leq \max I_{q} (h) \leq  - c_2/q .
\end{displaymath} 
Now, by the mean-value theorem, $I_{q} (h)$ has length $2h/ \vert B' (t_q)\vert$ for some $t_q \in I_{q} (h)$. But, using \eqref{goodbehaviour},  
we get:  
\begin{displaymath} 
B (t) \approx \frac{1}{t} \quad \text{ and} \quad  \vert B' (t) \vert \approx \frac{1}{t^2} \, \raise 1pt \hbox{,} 
\end{displaymath} 
so that $I_{q} (h)$ has length approximately $h t_{q}^2 \approx h/ q^2$ since $\vert t_q \vert \approx 1/ q$. Because of the modulus constraint, the only 
involved $q$'s are those for which $q \geq q_1$, where $q_1 \approx 2^{n/ \theta}$. Taking $n$ numerically large enough, we may assume that $q_1 > q_0$. 
We finally see that, for any $0 \leq j \leq 2^{n} - 1$, we have the lower bound: 
\begin{displaymath} 
m_{\varphi} (W_{n, j}) = m (E_{n, j}) 
\gtrsim \sum_{q \geq q_1} m (I_{q} (h)) 
\gtrsim \sum_{q \geq q_1} \frac{h}{q^2} 
\gtrsim \frac{h}{q_1} 
\gtrsim 2^{- n (1 + 1/ \theta)}.
\end{displaymath} 
It follows that: 
\begin{align*}
\sum_{n=1}^\infty \sum_{j=0}^{2^{n} - 1} \big[ 2^n m_{\varphi} (W_{n, j}) \big]^{p/2}  
& \gtrsim \sum_{n=1}^\infty \sum_{j=0}^{2^{n} - 1} \big[2^{n} 2^{- n (1 + 1/\theta)} \big]^{p/2} 
= \sum_{n=1}^\infty \sum_{j=0}^{2^{n} - 1} \big[2^{- n p/ 2\theta} \big] \\ 
& =\sum_{n=1}^\infty 2^{n (1 - p/ 2\theta)} 
= \infty , 
\end{align*}
if $p\leq 2 \theta$. Hence $C_{\varphi_2} \notin S_p$ for $p \leq 2 \theta$ by Proposition~\ref{lulu}. \par
\smallskip

A similar upper bound, and the membership of $C_{\varphi_2}$ in $S_p$ for $p > 2\theta$, would easily be proved along the same lines (and we will make use 
of that fact in Section~\ref{counterexamples}). But this will also follow from the more precise result on approximation numbers. To that effect, we shall borrow 
the following result from \cite{LIQUEROD}.
\begin{theoreme} [\cite{LIQUEROD}] \label{Auxiliary} 
Let $\varphi$ be a Schur function. Then the approximation numbers of  $C_{\varphi} \colon H^2 \to H^2$ have the upper bound:
\begin{equation} \label{Super} 
\qquad a_{n} (C_\varphi) \leq K \inf_{0 < h < 1} \bigg[(1 - h)^n + \sqrt{\frac{\rho_{\varphi} (h)}{h}} \bigg], \qquad n = 1, 2, \ldots \,. \quad 
\end{equation} 
\end{theoreme}

Applying this theorem to $\varphi_2$, which satisfies $\rho_{\varphi_2} (h) \leq K h^{1 + 1/\theta}$ as is clear from the preceding computations, would 
provide \emph{upper bounds} for $m_{\varphi} (W_{n, j})$ of the same order as the lower bounds obtained. Then choosing  
$h = H \log n / n$, where $H$ is a large constant ($H = 1/ 2 \theta$ will do) and using $1 - h \leq \e^{- h}$, we get from \eqref{Super}:
\begin{displaymath} 
a_{n} (C_{\varphi_2}) \leq K \bigg[n^{- H} + \bigg(\frac{\log n}{n}\bigg)^{1/ 2\theta} \bigg] 
\leq K \bigg(\frac{\log n}{n}\bigg)^{1/ 2\theta} . 
\end{displaymath} 

This ends the proof of Theorem~\ref{degrade}. \hfill $\square$
\bigskip \goodbreak

\noindent{\bf Remark:} Theorem~\ref{Auxiliary} of \cite{LIQUEROD} gives a very imprecise estimate on the approximation numbers of lens maps, as we 
noticed in that paper. On the other hand, when we apply it to a lens map spread by multiplication by the inner function $M$, we obtain an estimate which is close 
to being optimal, up to a logarithmic factor. This indicates that many phenomena have still to be understood concerning approximation numbers of 
composition operators.

 
\section{Lens maps as counterexamples} \label{counterexamples}

 Recall that an operator $T \colon X \to Y$ between Banach spaces is said to be Dunford-Pettis (in short DP) or completely continuous, if for any sequence 
$(x_n)$ which is weakly convergent to $0$, the sequence $(T x_n)$ is norm-convergent to $0$.  It is called weakly compact (in short $w$-compact) if the image 
$T(B_X)$ of the unit ball in $X$ is (relatively) weakly compact in $Y$. The identity map $i_1 \colon \ell_1 \to \ell_1$ is DP and not $w$-compact, by the Schur 
property of $\ell_1$ and its non-reflexivity. If $1 < p < \infty$, the identity map $i_p \colon \ell_p \to \ell_p$ is $w$-compact and not DP by the reflexivity of 
$\ell_p$ and the fact that the canonical basis $(e_n)$ of $\ell_p$ converges weakly to $0$, whereas $\Vert e_n\Vert_p = 1$.  Therefore, the two notions, clearly 
weaker than that of compactness, are not comparable in general. Moreover, when $X$ is reflexive, any operator $T \colon X\to Y$ is $w$-compact and any 
Dunford-Pettis operator $T \colon X \to Y$ is compact. \par 

Yet, in the context of composition operators $T = C_\varphi \colon X \to X$, with $X$ a non-reflexive Banach space of analytic functions, several results 
say that weak compactness of $C_\varphi$ implies its compactness. Let us quote some examples:
\par\smallskip

- $X = H^1$; this was proved by D. Sarason in 1990 (\cite{Sarason}); \par 
- $X = H^\infty$ and the disk algebra $X = A(\D)$ (A. \"Ulger \cite{Ulger} and R. Aron, P. Galindo and M. Lindstr\"om \cite{Aron}, independently; the 
first-named of us also gave another proof in \cite{Pascal-JOT}); \par 
- $X$ is the little Bloch space ${\mathscr B}_0$ (K. Madigan and A. Matheson \cite{M&M}); \par
- $X$ is the Hardy-Orlicz spaces $X = H^\psi$, when the Orlicz function $\psi$ grows more rapidly than power functions, namely when it satisfies the condition 
$\Delta^0$ (\cite{MEMOIRS}, Theorem~4.21, page~55); \par 
- $X = BMOA$ and $X = VMOA$ (J. Laitila, P. J. Nieminen, E. Saksman and H.-O. Tylli \cite{Finlande}).
\par\smallskip

Moreover, in some cases, $C_\varphi$ is compact whenever it is Dunford-Pettis (\cite{Pascal-JOT} for $X = H^\infty$ and \cite{MEMOIRS}, Theorem~4.21, 
page~55, for $X = H^\psi$, when the conjugate function of $\psi$ satisfies the condition $\Delta_2$). \par
\smallskip

The question naturally comes whether for \emph{any} non-reflexive Banach space $X$ of analytic functions on $\D$, \emph{every} weakly compact (resp. 
Dunford-Pettis) composition operator $C_\varphi \colon X \to X$ is actually compact. The forthcoming theorems show that the answer is \emph{negative} in 
general. Our spaces $X$ will be Hardy-Orlicz and Bergman-Orlicz spaces, so we first recall some definitions and facts about Orlicz spaces (\cite{MEMOIRS}, 
Chapters 2 and 3). 
\medskip

An \emph{Orlicz function} is a nondecreasing convex function $\psi \colon \R^+ \to \R^+$ such that $\psi (0) = 0$ and 
$\psi(\infty) = \infty$. Such a function is automatically continuous on $\R^+$. If $\psi (x) $ is not equivalent to an affine function, we must have 
$\psi (x) / x \mathop{\longrightarrow}\limits_{x \to \infty} \infty$. The Orlicz function $\psi$ is said to satisfy the $\Delta_2$-condition if 
$\psi (2x) / \psi (x)$ remains bounded. The conjugate function $\tilde \psi$ of an Orlicz function $\psi$ is the Orlicz function defined by:
\begin{displaymath} 
\tilde \psi (x) = \sup_{y \geq 0} \big(x y - \psi (y) \big).
\end{displaymath} 
For the conjugate function, one has the following characterization of $\Delta_2$  (see \cite{MEMOIRS}, page 7): $\tilde \psi$ has $\Delta_{2}$ if and only if , 
for some $\beta > 1$ and $x_0 > 0$, 
\begin{equation} \label{pareil} 
\qquad \quad \psi (\beta x) \geq 2 \beta \psi (x) \,, \qquad \text{for all } x \geq x_0.
\end{equation}
Let $(\Omega, \mathcal{A}, \P)$  be a probability space, and $L^0$ the space of measurable functions $f \colon \Omega \to \C$.  The Orlicz space 
$L^{\psi} = L^{\psi} (\Omega, \mathcal{A}, \P)$ is defined by
\begin{displaymath} 
L^{\psi} (\Omega, \mathcal{A}, \P) = \big\{f \in L^{0} \, ; \ \int_{\Omega} \psi (\vert f \vert / K) \, d \, \P <\infty \text{ for some } K > 0 \big\}.
\end{displaymath} 
This  is a Banach space for the Luxemburg norm:
\begin{displaymath} 
\Vert f \Vert_{\psi} = \inf  \big\{K > 0 \, ; \ \int_{\Omega} \psi (\vert f \vert / K) \, d \, \P \leq 1 \big\}.
\end{displaymath} 
The Morse-Transue space $ M^\psi$ (see \cite{MEMOIRS}, page~9) is the subspace of functions $f$ in $L^\psi$ for which 
$\int_{\Omega} \psi (\vert f \vert / K) \, d \, \P < \infty$ for every $K > 0$. It is the closure of $L^\infty$. One always has 
$(M^\psi)^* = L^{\tilde \psi}$ and $L^\psi = M^\psi$ if and only if $\psi$ has $\Delta_2$. When the conjugate function $\tilde \psi$ of $\psi$ 
has $\Delta_2$, the bidual of $M^\psi$ is then (isometrically isomorphic to) $L^\psi$. \par

Now, we can define the Hardy-Orlicz space $H^\psi$ attached to $\psi$ as follows. Take the probability space $(\T, \mathcal{B}, m)$ and recalling that 
$f_{r} (\e^{it}) = f (r \e^{it})$: 
\begin{displaymath} 
H^\psi = \{f \in {\cal H} (\D) \, ; \ \sup_{0 < r < 1} \Vert f_r \Vert_{L^{\psi} (m)} := \Vert f \Vert_{H^\psi} < \infty \}.
\end{displaymath} 
We refer to \cite{MEMOIRS} for more information on $H^\psi$. Similarly, we define (see \cite{MEMOIRS}) the Bergman-Orlicz space $B^\psi$, using this time 
the normalized area measure $A$, by:
\begin{displaymath} 
B^\psi = \{f \in {\cal H} (\D) \, ; \  \Vert f \Vert_{L^{\psi} (A)} := \Vert f \Vert_{L^\psi} < \infty \}.
\end{displaymath} 
If $\psi(x) = x^p$, $p \geq 1$, we get the usual Hardy and Bergman spaces $H^p$ and $B^p$. Those  spaces  are Banach spaces for any $\psi$, and Hilbert spaces 
for $\psi (x) = x^2$. The Hardy-Morse-Transue space $HM^\psi$ and Bergman-Morse-Transue space $BM^\psi$ are defined by 
$HM^\psi = H^\psi \cap M^\psi$ and $BM^\psi = B^\psi \cap M^\psi$. When the conjugate function of $\psi$ has $\Delta_2$, the bidual of $HM^\psi$ is 
(isometrically isomorphic to) $H^\psi$ (\cite{MEMOIRS}, page~10).
\par\bigskip

We can now state our first theorem.
\begin{theoreme} \label{Counterexample one} 
There exists a Schur function $\varphi$ and an Orlicz function $\psi$ such that $H^\psi$ is not reflexive and the composition operator 
$C_\varphi \colon H^\psi \to H^\psi$ is weakly-compact and Dunford-Pettis, but is not compact.
\end{theoreme}

\noindent{\bf Proof.} First take for $\varphi$ the lens map $\varphi_{1/ 2}$ which in view of (\ref{Edensor case}) of Lemma~\ref{simplet} satisfies, for 
some constant $K > 1$: 
 \begin{equation} \label{Goodminoration} 
\qquad \quad \rho_\varphi (h) \geq K^{-1} h^2, \qquad 0 < h < 1.
\end{equation}
We now recall the construction of an Orlicz function made in \cite{STUDIA}. Let  $(x_n)$ be a the sequence of positive numbers defined as follows: 
$x_1 = 4$ and then, for every integer $n \geq 1$, $x_{n+1} > 2x_n$ is the abscissa of the second intersection point of the parabola $y = x^2$ with the straight line 
containing $(x_n,  x_{n}^2)$ and $(2 x_n, x_{n}^4)$; equivalently $x_{n+1} = x_{n}^3 - 2 x_n$. We now define our Orlicz function $\psi$ by 
$\psi (x) = 4 x$ for $0 \leq x \leq 4$ and, for $n \geq 1$, by: 
\begin{equation}\label{induction} 
\begin{split}
\psi (x_n) & = x_{n}^2 , \\    
\psi \text{ affine\ between } x_n \text{ and } & x_{n+1}, \text{ so that } \psi (2 x_n) = x_{n}^4. 
\end{split}
\end{equation} 
Observe that $\psi$ does not satisfy the $\Delta_2$-condition, since $\psi (2x_n) = [\psi (x_n)]^2$. It clearly satisfies (since $\psi^{- 1}$ is concave): 
\begin{equation} \label{inbetween} 
\begin{split}
& \qquad x^2 \leq \psi (x) \leq x^4 \qquad \quad \text{ for } x \geq 4, \\ 
& \qquad \psi^{- 1} (K x) \leq K \psi^{- 1} (x) \quad \text{ for any } x > 0,  K > 1.
\end{split}
\end{equation}
Therefore, it has a moderate growth, but a highly irregular behaviour, which will imply the results we have in view. Indeed, let $y_n = \psi (x_n)$ and 
$h_n = 1 / y_n$. We see from \eqref{Goodminoration}, \eqref{induction} and \eqref{inbetween} that:
\begin{equation} \label{criterus} 
D (h_n) \mathop{=}^{def} \frac{\psi^{- 1}(1 / h_n)}{\psi^{- 1}(1 / \rho_{\varphi} (h_n))} 
\geq \frac{\psi^{- 1}(1 / h_n)}{\psi^{- 1}(K / h_{n}^2)} 
= \frac{\psi^{- 1} (y_n)}{\psi^{- 1} (K y_{n}^2)} 
\geq \frac{x_n}{2 K x_n} = \frac{1}{2 K}. 
\end{equation}
Thus, we have $\limsup_{h \to 0^+} D (h) > 0$. By \cite{MEMOIRS}, Theorem~4.11 (see also \cite{MathAnn}, comment before Theorem~5.2), 
$C_\varphi$ is not compact. \par
\smallskip

On the other hand, let $j_{\psi, 2} \colon H^\psi \to H^2$ and $j_{4, \psi} \colon H^4 \to H^\psi$ be the natural injections, which are continuous,  
thanks to \eqref{inbetween}.  We have the following diagram:
\begin{displaymath} 
H^\psi \mathop{\longrightarrow}^{j_{\psi, 2}} H^2 \mathop{\longrightarrow}^{C_\varphi} H^4  \mathop{\longrightarrow}^{j_{4, \psi}} H^\psi . 
\end{displaymath} 
The second map is continuous as a consequence of \eqref{Edensor case} and of a result of P. Duren (\cite{Duren}; see also \cite{DUREN-livre}, Theorem~9.4, 
page~163),  which extends Carleson's embedding theorem (see also \cite{MEMOIRS}, Theorem~4.18). Hence 
$C_\varphi = j_{4, \psi} \circ C_\varphi \circ j_{\psi, 2}$ factorizes through a reflexive space ($H^2$ or $H^4$) and is therefore $w$-compact. 
\par \smallskip

To prove that $C_\varphi$ is Dunford-Pettis, we use the following result of \cite{TRANS} (Theorem~2.1):
\begin {theoreme} [\cite{TRANS}] \label{use} 
Let $\varphi$ be a Schur function and $\Phi$ be an Orlicz function. Assume that, for some $A > 0$, one has: 
\begin{equation} \label{gotobed} 
\qquad \quad \sup_{0 < t \leq h} \frac{\rho_{\varphi} (t)}{t^2} \leq \frac{1 / h^2}{\Phi \big(A \Phi^{-1}(1 / h^2) \big)} \, \raise 1pt \hbox{,} 
\qquad 0 < h < 1. 
\end{equation} 
Then, the canonical inclusion $j_{\Phi, \varphi} \colon B^\Phi \to L^{\Phi} (m_\varphi)$ is continuous. \par
In particular, it is continuous for any Orlicz function $\Phi$ if $\rho_{\varphi} (h) = O\, (h^2)$.
\end{theoreme} 

Now, let $J_\psi \colon H^\psi \to B^\psi$ be the canonical inclusion, and consider the following diagram:
\begin{displaymath} 
H^\psi \mathop{\longrightarrow}^{J_\psi} B^\psi \mathop{\longrightarrow}^{j_{\psi, \varphi}} L^{\psi} (m_\varphi).
\end{displaymath} 
The first map is Dunford-Pettis, by \cite{STUDIA}, Theorem~4.1. The second map is continuous by \eqref{Edensor case} and \eqref{gotobed}.  
Clearly, being Dunford-Pettis is an ideal property (if either $u$ or $v$ is Dunford-Pettis, so is $vu$). Therefore, $j_{\psi, \varphi} \circ J_\psi$ is 
Dunford-Pettis, and this amounts to say that $C_\varphi \colon H^\psi \to H^\psi$ is Dunford-Pettis. 
\par\smallskip

Now,  the non-reflexivity of $H^\psi$ follows automatically, since $C_\varphi$ is  Dunford-Pettis but not compact. \par \smallskip  

This ends the proof of Theorem~\ref{Counterexample one}.  \hfill $\square$
\bigskip
 
Theorem~\ref{Counterexample one} admits the following variant.  

\begin{theoreme} \label{Counterexample two} 
There exist a Schur function $\varphi$ and an Orlicz function $\chi$ such that $H^\chi$  is not reflexive and the composition operator 
$C_\varphi \colon H^\chi \to H^\chi$ is weakly compact and not Dunford-Pettis; in particular it is not compact.
\end{theoreme}

\noindent{\bf  Proof .} We use the same Schur function $\varphi = \varphi_{1 / 2}$, but we replace $\psi$ by the function $\chi$ defined by 
$\chi (x) = \psi (x^2)$. Let $A > 1$. Observe that, in view of \eqref{inbetween}, 
\begin{displaymath} 
 \frac{\chi (A x)}{[\chi (x)]^2} = \frac{\psi (A^2 x^2)}{[\psi (x^2)]^2} \leq \frac{A^8 x^8}{x^8}=A^8.
\end{displaymath} 
By \cite{STUDIA}, Proposition~4.4, $J_{\chi} \colon H^{\chi} \to B^{\chi}$ is $w$-compact, and  we can see 
$C_\varphi \colon H^\chi \to H^\chi$ as the canonical inclusion $j \colon H^\chi \to L^{\chi} (m_\varphi)$. Hence Theorem~\ref{use} and the diagram: 
\begin{displaymath} 
j = j_{\chi, \varphi} \circ J_\chi \colon 
H^\chi  \mathop{\longrightarrow}^{J_\chi} B^\chi \mathop{\longrightarrow}^{j_{\chi, \varphi}} L^{\chi} (m_\varphi)  
\end{displaymath} 
show that $C_\varphi \colon H^\chi \to H^\chi$ is $w$-compact as well. \par
\smallskip

Now, to prove that $C_\varphi$ is not Dunford-Pettis, we cannot use \cite{STUDIA}, as in the proof of Theorem~\ref{Counterexample one}, but we follow 
the lines of Proposition~3.1 of \cite{STUDIA}. Remark first that, by definition, the function $\chi$ satisfies, for $\beta = 2$,  the following inequality:
\begin{displaymath} 
\chi (\beta x) = \psi (4 x^2) \geq 4 \psi (x^2) = 2 \beta \chi (x) ;
\end{displaymath} 
hence, by \eqref{pareil}, this implies that the conjugate function of $\chi$ verifies the $\Delta_2$-condition. \par

Let $x_n$ be as in \eqref{induction}, and set:
\begin{displaymath} 
u_n = \sqrt {x_n} \qquad \text{and} \quad  A = \sqrt 2
\end{displaymath} 
so that 
\begin{equation} \label{tantum} 
\chi(Au_n)=\big[\chi( u_n)\big]^2=x_{n}^4. 
\end{equation}
Finally, let:
\begin{displaymath} 
r_n = 1 - \frac{1}{\chi ( u_n)} \qquad \text{and} \quad f_{n} (z) = u_n \bigg(\frac{1 - r_n}{1 - r_{n} z} \bigg)^2.
\end{displaymath} 
By (\cite{MEMOIRS}, Corollary~3.10), $\Vert f_n \Vert _{H^\chi} \leq 1$ and $f_n$ tends to $0$ uniformly on compact subsets of $\D$; that implies that 
$f_n \to 0$ weakly in $H^\chi$ since the conjugate function of $\chi$ has $\Delta_2$ (\cite{MEMOIRS}, Proposition~3.7). On the other hand, if 
$K_n = \Vert f_n \Vert_{L^{\chi} (m_\varphi)}$, mimicking the computation of (\cite{STUDIA}, Proposition 3.1), we get:
\begin{equation} \label{imite} 
1 = \int_{\D} \chi (\vert f_n \vert / K_n) \, dm_\varphi 
\geq (1 - r_n)^2 \chi (\alpha u_n / 4 K_n) 
\end{equation} 
for some $0 < \alpha < 1$ independent of $n$, where we used the convexity of $\chi$ and the fact that the lens map $\varphi$ satisfies, 
by\eqref{Goodminoration}: 
\begin{displaymath} 
m_{\varphi} (\{z \in\D \, ; \ \vert 1 - z \vert \leq 1 - r_n \}) \geq \alpha (1 - r_n)^2.
\end{displaymath} 
In  view of \eqref{tantum}, \eqref{imite} reads as well: 
\begin{displaymath} 
\chi (\alpha u_n / 4 K_n) \leq \chi^{2} ( u_n) = \chi (A u_n), 
\end{displaymath} 
so that:
\begin{equation} \label{ergo} 
\Vert j (f_n) \Vert_{L^{\chi} (m_\varphi)} = K_n \geq \alpha / 4 A. 
\end{equation}
This shows that  $j \colon H^\chi \to L^{\chi} (m_\varphi)$ and therefore also $C_\varphi \colon H^\chi \to H^\chi$ are  not Dunford-Pettis. 
\par\smallskip

It remains to show that $H^\chi$ is not reflexive. We shall prove below a more general result, but here, the conjugate function $\tilde \chi$ of $\chi$ 
satisfies the $\Delta_2$ condition, as we saw. Hence $H^\chi$ is the bidual of $HM^\chi$. Since $\chi$ fails to satisfy the $\Delta_2$-condition, we know 
that $L^\chi \neq M^\chi$. Let $u \in L^\chi \setminus M^\chi$,with  $u\geq 1$. Let $f$ be the associated outer function, namely: 
\begin{displaymath} 
f (z) = \exp\Big(\frac{1}{2\pi} \int_{0}^{2\pi} \frac{\e^{it} + z}{\e^{it} - z} \log u (t) \, dt\Big).
\end{displaymath} 
One has $\vert f^*\vert  = u$ almost everywhere, with the notations of \eqref{radial}, and hence $f \in H^\chi \setminus HM^\chi$. It follows that  
$H^\chi\neq HM^\chi$. Hence $HM^\chi$ is not reflexive, and therefore $H^\chi$ is not reflexive either. \hfill $\square$
\bigskip

As promised, we give the general result on non-reflexivity.
\begin{proposition}
Let $\psi$ be an Orlicz function which does not satisfy the $\Delta_2$-condition. Then neither $H^\psi$ nor $B^\psi$ is reflexive.
\end{proposition}

\noindent{\bf Proof.} We only give the proof for $B^\psi$ because it is the same for $H^\psi$. \par
Since $\psi$ does not satisfy $\Delta_2$ there is a sequence $(x_n)$ of positive numbers, tending to infinity, such that $\psi (2 x_n) / \psi (x_n)$ tends to 
infinity. Let $r_n \in (0, 1)$ such that $(1 - r_n)^2 = 1/ \psi (2 x_n)$ and set:
\begin{displaymath} 
q_n (z) = \frac{(1 - r_n)^4}{(1 - r_n z)^4} \cdot
\end{displaymath} 
One has $\|q_n \|_\infty = 1$ and $\| q_n \|_1 = \frac{(1 - r_n)^2}{(1 + r_n)^2} \leq (1 - r_n)^2$. On the other hand, on the pseudo-Carleson 
window $S (1, 1 - r_n)$, one has $|1 - r_n z| \leq (1 - r_n) + r_n |1 - z| \leq (1 - r_n) + r_n (1 - r_n) = 1 - r_n^2 \leq 2 (1 - r_n)$; 
hence $|q_n (z) | \geq 1/16$. It follows that:
\begin{align*}
1 
& = \int_\D \psi \bigg(\frac{|q_n|}{\|q_n\|_\psi} \bigg)\, dA 
\geq \int_{S (1, 1 - r_n)} \psi \bigg( \frac{|q_n|}{\|q_n\|_\psi} \bigg) \, dA \\
& \geq A [S (1, 1 - r_n)]\, \psi \bigg( \frac{1}{16 \, \|q_n\|_\psi} \bigg) 
\geq \frac{1}{3} \, (1 - r_n)^2 \psi \bigg( \frac{1}{16 \, \|q_n\|_\psi} \bigg) \\
& \geq (1 - r_n)^2 \psi \bigg( \frac{1}{48 \, \|q_n\|_\psi} \bigg) 
= \frac{1}{\psi (2 x_n)} \, \psi \bigg( \frac{1}{48 \, \|q_n\|_\psi} \bigg) ; 
\end{align*} 
hence $\psi (1/ [48 \, \|q_n\|_\psi]) \leq \psi (2 x_n)$, so $1/ (48 \, \|q_n\|_\psi) \leq 2 x_n$ and $96 \, x_n \, \|q_n\|_\psi \geq 1$. 
\par\smallskip

Set now $f_n = q_n / \|q_n\|_\psi$; one has $\|f_n\|_\psi = 1$ and (using that $\psi ( x_n \, |q_n (z)| ) \leq |q_n (z)|\, \psi (x_n)$, by convexity, since 
$|q_n (z)| \leq 1$):
\begin{align*} 
\int_\D \psi \bigg( \frac{|f_n|}{96} \bigg) \, dA 
& = \int_\D \psi \bigg( \frac{x_n \, |q_n|}{96 \, x_n \, \|q_n\|_\psi} \bigg) \, dA 
\leq \int_\D \psi (x_n \, |q_n|) \, dA \\ 
& \leq \psi (x_n) \int_\D |q_n| \, dA \,, \\ 
& \leq \psi (x_n) \, (1 - r_n)^2 = \frac{\psi (x_n)}{\psi (2 x_n)} \mathop{\longrightarrow}_{n \to \infty} 0. 
\end{align*}  
By \cite{critere}, Lemma~11, that implies that the sequence $(f_n)$ has a subsequence equivalent to the canonical basis of $c_0$ 
and hence $B^\psi$ is not reflexive. \qed 
\bigskip 

We finish by giving a counterexample using Bergman-Orlicz spaces instead of Hardy-Orlicz spaces.
\begin{theoreme} \label{Counterexample three} 
There exists a Schur function $\varphi$ and an Orlicz function $\psi$ such that the space $B^\psi$ is not reflexive and 
the composition operator $C_\varphi \colon B^\psi \to B^\psi$ is  weakly-compact but not compact. 
\end{theoreme}
\noindent{\bf Proof.} We use again the Orlicz function $\psi$ defined by \eqref{induction} and the Schur function $\varphi = \varphi_{1/2}$. 
The space $B^\psi$ is not reflexive since $\psi$ does not satisfy the condition $\Delta_2$. \par
We now need an estimate similar to \eqref{Edensor case} for $\varphi_\theta$, namely: 
\begin{equation} \label{Stefan case} 
\rho_{\varphi, 2} (h) := \sup_{|\xi| = 1} A [\{z \in \D \, ; \ \varphi (z) \in S (\xi, h)\}] \approx h^{2 / \theta} . 
\end{equation}
The proof of \eqref{Stefan case} is  best seen by passing to the right half-plane with the measure $A_{\gamma_\theta}$ which is locally equivalent 
to the Lebesgue planar measure $A$; we get 
$\rho_{\varphi, 2} (h) \geq A (\{\vert z\vert^\theta \leq h\} \cap \H) \geq K h^{2 / \theta}$ and the upper bound in \eqref{Stefan case} is proved 
similarly. \par

We now see that $C_\varphi \colon B^\psi \to B^\psi$ is not compact as follows. We use the same $x_n$ as in \eqref{induction} and set 
$y_n = \psi (x_n)$, $k_n = 1 / \sqrt{y_n}$. We notice that, since $ \rho_{\varphi, 2} (h) \geq K^{-1} h^4$ (with $K > 1$) in view of 
\eqref{Stefan case}, we have: 
\begin{displaymath} 
E(k_n) 
\mathop{=}^{def} \frac{\psi^{-1} (1 / k_{n}^2 )}{\psi^{-1} (1 / \rho_{\varphi, 2} (k_n))}  
\geq \frac{\psi^{-1}(1 / k_{n}^2)}{\psi^{-1}(K / k_{n}^4)} 
= \frac{\psi^{-1}(y_n)}{\psi^{-1}(K y_{n}^2)} 
\geq \frac{x_n}{2 K x_n} = \frac{1}{2 K} \, \raise 1pt \hbox{,} 
\end{displaymath} 
so that  
\begin{displaymath} 
\limsup_{k \to 0^+} E (k) > 0, 
\end{displaymath} 
and this implies that $C_\varphi \colon B^\psi \to B^\psi$ is not compact (\cite{TRANS}, Theorem~3.2). To see that 
$C_\varphi \colon B^\psi\to B^\psi$ is $w$-compact, we use the diagram: 
\begin{displaymath} 
B^\psi \mathop{\longrightarrow}^{j_{\psi, 2}} B^2 \mathop{\longrightarrow}^{C_\varphi} B^4 \mathop{\longrightarrow}^{j_{4, \psi}} B^\psi  
\end{displaymath} 
as well as \eqref{Stefan case}, which gives $\rho_{\varphi, 2} (h) \leq K h^4$. A result of W. Hastings (\cite{HAS}) now implies the continuity of 
the second map. This diagram shows that $C_\varphi$ factors through a reflexive space ($B^2$ or $B^4$), and is therefore $w$-compact.
 \hfill $\square$


\bigskip

\vbox{\small \noindent{\it
{\rm Pascal Lef\`evre}, Univ Lille Nord de France, \\
U-Artois, Laboratoire de Math\'ematiques de Lens EA~2462, \\
F\'ed\'eration CNRS Nord-Pas-de-Calais FR~2956, \\
Facult\'e des Sciences Jean Perrin,\\
Rue Jean Souvraz, S.P.\kern 1mm 18, \\ 
F-62\kern 1mm 300 LENS, FRANCE \\ 
pascal.lefevre@euler.univ-artois.fr
\smallskip

\noindent
{\rm Daniel Li}, Univ Lille Nord de France, \\
U-Artois, Laboratoire de Math\'ematiques de Lens EA~2462, \\
F\'ed\'eration CNRS Nord-Pas-de-Calais FR~2956, \\
Facult\'e des Sciences Jean Perrin,\\
Rue Jean Souvraz, S.P.\kern 1mm 18, \\ 
F-62\kern 1mm 300 LENS, FRANCE \\ 
daniel.li@euler.univ-artois.fr
\smallskip

\noindent
{\rm Herv\'e Queff\'elec}, Univ Lille Nord de France, \\
USTL, Laboratoire Paul Painlev\'e U.M.R. CNRS 8524, \\ 
F\'ed\'eration CNRS Nord-Pas-de-Calais FR~2956, \\
F-59\kern 1mm 655 VILLENEUVE D'ASCQ Cedex, 
FRANCE \\ 
Herve.Queffelec@univ-lille1.fr
\smallskip

\noindent
{\rm Luis Rodr{\'\i}guez-Piazza}, Universidad de Sevilla, \\
Facultad de Matem\'aticas, Departamento de An\'alisis Matem\'atico {\rm \&} IMUS,\\ 
Apartado de Correos 1160,\\
41\kern 1mm 080 SEVILLA, SPAIN \\ 
piazza@us.es\par}
}


\begin{thebibliography}{99}

\bibitem{Aron} R. Aron, P. Galindo, M. Lindstr\"om, 
Compact homomorphisms between algebras of analytic functions, 
Studia Math. 123 (1997), no. 3, 235--247.

\bibitem{CA-ST-livre} B. Carl, M. Stephani, 
Entropy, Compactness and the Approximation of Operators, 
Cambridge University Press (1990). 

\bibitem{Duren} P. Duren, 
Extension of a theorem of Carleson, 
Bull. Amer. Math. Soc. 75 1969 143--146.

\bibitem{DUREN-livre} P. Duren,
Theory of $H^p$-spaces, Second edition, 
Dover Publications (2000).

\bibitem{HAS} W. H. Hastings, 
A Carleson measure theorem for Bergman spaces, 
Proc. Amer. Math.Soc. 52 (1975), 237--241. 

\bibitem{Finlande} J. Laitila, P. J. Nieminen, E. Saksman, H.-O. Tylli, 
Compact and Weakly Compact Composition Operators on BMOA, 
\emph{to appear} in Complex Anal. Operator Theory {\small DOI: 10.1007/s11785-011-0130-9}.

\bibitem{Pascal-JOT} P. Lef\`evre, 
Some characterizations of weakly compact operators in $H^\infty$ and on the disk algebra. Application to composition operators, 
J. Operator Theory 54 (2005), no. 2, 229--238.

\bibitem{critere} P. Lef\`evre, D. Li, H. Queff\'elec, L. Rodr{\'\i}guez-Piazza, 
A criterion of weak compactness for operators on subspaces of Orlicz spaces,  
Journal of Function Spaces and Applications 6, No.~3 (2008), 277--292.

\bibitem{JFA} P. Lef\`evre, D. Li, H. Queff\'elec, L. Rodr{\'\i}guez-Piazza, 
Some examples of compact composition operators on $H^2$, 
J. Funct. Anal. 255, No.11 (2008), 3098--3124. 

\bibitem{JMAA} P. Lef\`evre, D. Li, H. Queff\'elec, L. Rodr{\'\i}guez-Piazza, 
Compact composition operators on $H^2$ and Hardy-Orlicz spaces, 
J. Math. Anal. Appl. 354 (2009), 360--371. 

\bibitem{MEMOIRS} P. Lef\`evre, D. Li, H. Queff\'elec, L. Rodr{\'\i}guez-Piazza, 
Composition operators on Hardy-Orlicz spaces,  
Memoirs Amer. Math. Soc. 207 (2010), no.~974. 

\bibitem{MathAnn} P. Lef\`evre, D. Li, H. Queff\'elec, L. Rodr{\'\i}guez-Piazza, 
Nevanlinna counting function and Carleson function of analytic maps, 
Math. Ann.  351 (2011), 305--326. 

\bibitem{STUDIA} P. Lef\`evre, D. Li, H. Queff\'elec, L. Rodr{\'\i}guez-Piazza, 
The canonical injection of the Hardy-Orlicz  space $H^\psi$ into the Bergman-Orlicz space $\mathcal{B}^\psi$, 
Studia Math. 202 (1) (2011), 123--144.

\bibitem{TRANS} P. Lef\`evre, D. Li, H. Queff\'elec, L. Rodr{\'\i}guez-Piazza, 
Compact composition operators on Bergman-Orlicz spaces, 
\emph{preprint} {\small arXiv : 0910.5368}. 

\bibitem{LIQUEROD}  D. Li, H. Queff\'elec, L. Rodr{\'\i}guez-Piazza, 
On approximation numbers of composition operators.   

\bibitem{LUECKING} D. H. Luecking, 
Trace ideal criteria for Toeplitz operators, 
J. Funct. Anal. 73 (1987), 345--368.


\bibitem{M&M} K. Madigan, A. Matheson, 
Compact composition operators on the Bloch space, 
Trans. Amer. Math. Soc. 347 (1995), no. 7, 2679--2687.

\bibitem{Pis-livre} G. Pisier, 
The volume of convex bodies and Banach space geometry, 
Cambridge University Press (1989).

\bibitem {Parf} O. G. Parfenov, 
Estimates of the singular numbers of the Carleson embedding operator, 
Math. USSR Sbornik 59 (2) (1988), 497--511. 

\bibitem{Sarason} D. Sarason, 
Weak compactness of holomorphic composition operators on $H^1$, 
Functional analysis and operator theory (New Delhi, 1990), 75--79, Lecture Notes in Math., 1511, Springer, Berlin (1992).

\bibitem{Shap-livre} J. H. Shapiro, 
Composition Operators and Classical Function Theory, 
Universitext, Tracts in Mathematics, Springer-Verlag, New York (1993). 

\bibitem{SHTA} J. H. Shapiro, P. D. Taylor, 
Compact, nuclear, and Hilbert-Schmidt composition operators on $H^2$, 
Indiana Univ. Math. J. 23 (1973), 471--496. 

\bibitem{Ulger} A. \"Ulger, 
Some results about the spectrum of commutative Banach algebras under the weak topology and applications, 
Monatsh. Math. 121 (1996), no.~4, 353--379.

\bibitem{Zhu-livre} K. Zhu, 
Operator Theory in Function Spaces, 
Mathematical Surveys and Monographs, American Math. Society, Vol. 138,  New York (2007). 

\end{thebibliography}
\end{document}